\documentclass[12pt,a4paper,british,fleqn]{scrartcl}
\usepackage[T1]{fontenc}
\usepackage[latin9]{inputenc}
\usepackage{color}
\usepackage{babel}
\usepackage{amsthm}
\usepackage{amsmath}
\usepackage{amssymb}
\usepackage[unicode=true,pdfusetitle,
 bookmarks=true,bookmarksnumbered=false,bookmarksopen=false,
 breaklinks=true,pdfborder={0 0 1},backref=false,colorlinks=true]
 {hyperref}

\makeatletter

\pdfpageheight\paperheight
\pdfpagewidth\paperwidth

\numberwithin{equation}{section}
\numberwithin{figure}{section}
\newcommand{\lyxaddress}[1]{
\par {\raggedright #1
\vspace{1.4em}
\noindent\par}
}

\@ifundefined{date}{}{\date{}}
\makeatother

\begin{document}

\title{Comment on \emph{``The Stochastic Nonlinear Schrödinger Equation
in $H^{1}$''}}

\author{Torquil Macdonald Sørensen}

\maketitle

\lyxaddress{\begin{center}
Centre of Mathematics for Applications\\
University of Oslo\\
NO-0316 Oslo, Norway\\
Email: t.m.sorensen@cma.uio.no, torquil@gmail.com
\par\end{center}}
\begin{abstract}
The paper \emph{``The Stochastic Nonlinear Schrödinger Equation in
$H^{1}$''} \cite{debouard2003} gives an existence proof for a stochastic
nonlinear Schrödinger equation with multiplicative noise. We point
out two mistakes that draw the validity of the proof into question.\\
\\
Keywords (MSC2010): 35Q41 Time-dependent Schrödinger equations, Dirac
equations; 35R60 Partial differential equations with randomness, stochastic
partial differential equations; 35G20 Nonlinear higher-order equations.
\end{abstract}

\section{Regarding the proof of \cite[Theorem 4.1]{debouard2003}\label{sec:theorem_4_1}}

Consider \cite[Theorem 4.1]{debouard2003}, concerning solution existence
for the following stochastic nonlinear Schrödinger equation with multiplicative
noise,
\begin{equation}
idu-(\Delta u+\lambda|u|^{2\sigma}u)dt=udW-\frac{i}{2}uF_{\phi}dt\,,\label{eq:spde}
\end{equation}
describing a stochastic process $u$ on $\mathbb{R}^{n}$. We refer
to \cite{debouard2003} for additional information about the mathematical
details. Here we only describe the bare minimum of details that are
necessary to describe our objection to the proof that is provided.

In the theorem, the following $n$-dependent parameter ranges for
$\sigma$ are assumed,
\begin{equation}
\begin{cases}
0<\sigma & ,n=1,2\\
0<\sigma<2 & ,n=3\\
\frac{1}{2}\leq\sigma<\frac{2}{n-2}\quad\mbox{or}\quad\sigma<\frac{1}{n-1} & ,n\geq4\,.
\end{cases}\label{eq:sigma_range}
\end{equation}
The theorem essentially states that an \emph{admissible pair} of Lebesgue
space exponents $(r,p)$ exists such that \eqref{eq:spde} has a unique
solution in a certain function space characterised by $(r,p)$. Admissibility
for $(r,p)$ is defined as
\begin{equation}
r\geq2,\quad\frac{2}{r}=n\left(\frac{1}{2}-\frac{1}{p}\right)\,.\label{eq:admissible_pair}
\end{equation}
Dual Lebesgue space exponents are denoted by primed quantities, and
are defined by the equation
\[
\frac{1}{p}+\frac{1}{p'}=1\,.
\]

The proof given is for the special case $\sigma\geq1/2$. In the proof,
a second admissible pair $(\gamma,s)$ is introduced after \cite[Equation (4.17)]{debouard2003}.
The parameters $s'$ and $p$ are related through another parameter
$q$ which arises in the proof,
\begin{equation}
\frac{1}{s'}=\frac{2\sigma}{q}+\frac{1}{p}\,,\label{eq:s'_q_p}
\end{equation}
which is described prior to \cite[Equation (4.18)]{debouard2003}.
The parameter $q$ arises due to the use of the Sobolev embedding
$H^{1}(\mathbb{R}^{n})\subset L^{q}(\mathbb{R}^{n})$. It is claimed
that the embedding holds because ``$q<2n/(n-3)<2n/(n-2)$''. However,
the second part of this inequality is incorrect, and therefore the
Sobolev embedding is used without proper justification.

\section{Regarding the proof of \cite[Lemma 4.3]{debouard2003}\label{sec:lemma_4_3}}

In the proof of \cite[Lemma 4.3]{debouard2003}, in the second estimate
on page 121, the interpolation inequality for $L^{p}$-spaces, followed
by the Sobolev embeddings $H^{1}(\mathbb{R}^{n}),W^{1,p}(\mathbb{R}^{n})\subset L^{q}(\mathbb{R}^{n})$
were used, where the parameter $q$ was ``as above''. Therefore,
for the same reason, the Sobolev embeddings used here also lack a
proper justification.

\section{Conclusions}

We believe that the errors described here put into question the validity
of the existence proof provided in \cite{debouard2003} for the SPDE
in the multiplicative case.

\bibliographystyle{utphys}
\bibliography{references}

\end{document}